\documentclass{amsart}
\usepackage[utf8]{inputenc}

\title{Motivic Splittings For Symmetric Matrices}
\author{Anubhav Nanavaty}

\thanks{This research was supported in part by NSF Grant No. DMS-1944862.}
\usepackage{spectralsequences}
\usepackage{amsmath}
\usepackage{graphicx}
\usepackage{amssymb}
\usepackage{amsthm}
\usepackage{mathtools}
\usepackage{adjustbox}
\usepackage{dsfont}
\usepackage{MnSymbol}
\usepackage{mathrsfs}
\usepackage{stmaryrd}
\usepackage{xcolor}
\usepackage[all]{xy}
\usepackage[mathcal]{eucal}
\usepackage{verbatim}  %%includes comment environment
\usepackage{fullpage}  %%smaller margins
\usepackage{hyperref}
\usepackage{setspace}
\usepackage{tikz-cd}
\usetikzlibrary{positioning}
\newcommand{\A}{\mathbb{A}}

\newcommand{\C}{\mathbb{C}}

\newcommand{\F}{\mathbb{F}}
\newcommand{\G}{\mathbb{G}}

\renewcommand{\L}{\mathbb{L}}

\newcommand{\N}{\mathbb{N}}
\renewcommand{\P}{\mathbb{P}}
\newcommand{\Q}{\mathbb{Q}}

\newcommand{\Z}{\mathbb{Z}}

\renewcommand{\phi}{\varphi}

\newcommand{\Sym}{\mathrm{Sym}}
\usepackage{rotating}
\usepackage{amssymb}

\newtheorem{theorem}{Theorem}[section]
\newtheorem{proposition}[theorem]{Proposition}
\newtheorem{lemma}[theorem]{Lemma}
\newtheorem{corollary}[theorem]{Corollary}

\newtheorem{remark}[theorem]{Remark}
\theoremstyle{definition}
\newtheorem{definition}[theorem]{Definition}

\numberwithin{equation}{subsection}
\begin{document}
\maketitle
\begin{abstract}
We show that the space of symmetric matrices of a fixed rank $k$ over a field $K$ of characteristic not equal to 2 is split Tate. We do this by promoting the point-counting strategy of MacWilliams over finite fields to a filtration of the locus of rank $\leq k$ symmetric matrices that is independent of the field. This filtration immediately allows for a computation of their isomorphism classes in the Grothendieck ring of varieties in terms of the Lefschetz motive. We then promote this computation to prove that the space of symmetric matrices of a fixed rank $k$ are split Tate in Voevodsky's category of motives in characteristic $0$ and Kelly's category of motives in characteristic $p$.
\end{abstract}
\section{Introduction}
Fix a field $K$ of characteristic not equal to $2$, and let $\text{Var}_K$ denote the category of $K$-varieties, or separated schemes of finite type over $K$. Let $\Sym^{n,k}$ (resp. $\Sym^{n,[k,\ell]}$) denote the variety of $n\times n$ symmetric matrices of rank $k$ (resp. in-between $k$ and $\ell$) over $K$. The original goal of this note was to understand the following properties of these varieties:
\begin{enumerate}
    \item If $K=\F_q$, the \'{e}tale cohomology with coefficients in $\Q_\ell$, for $\ell\neq q$
    \item If $K\hookrightarrow\C$, the mixed Hodge structure
\end{enumerate}
However, it became clear that many of the computations were actually being done in Voevodsky's category of motives $DM^{eff}_{gm,Nis}(K,\F)$ for any field $\F$ if $char(K)=0$ as defined in \cite{Vo1}, or Kelly's category of motives $DM^{eff}_{gm,\ell dh}(K,\F)$ if $char(K)=p$ and $p\in \F$ is invertible as defined in \cite{K}. We let $C^{gm}$ denote $DM^{eff}_{gm,\ell dh}(K,\F)$ or $DM^{eff}_{gm,Nis}(K,\F)$ depending on the context. We show the following:
\begin{theorem}\label{THE THEOREM}
    Let $K$ be a field of characteristic not equal to $2$. In $C^{gm}$, for all $n,k,\ell \geq 0$, the Voevodsky motives $M_{gm}(\Sym^{n, k})$ and $M_{gm}^c(\Sym^{n, [k,\ell]})$ are \emph{split Tate}. In other words, they are direct sums of the Tate motives $\F(i)[j]$, for $i,j\in\Z$.
\end{theorem}
Being split Tate is an extremely restrictive property, as the motives in $C^{gm}$ are often defined by non-trivial extensions. For example, the morphisms between two Tate motives are given by rational pieces of algebraic $K$ theory of the underlying field (see 6.2-6.3 of \cite{ClementDupont} for exposition). As the motives in \ref{THE THEOREM} are split Tate, results for the cohomology of these motives after passing to singular cohomology or \'{e}tale cohomology will provide formulas for the splitting in $C^{gm}$. We believe that these splitting formulas will be sensitive to the ground field, as various Chern classes that appear in the proof of \ref{THE THEOREM} may vanish due to the particular characteristic. 

There are many known instances in which the Voevodsky motive of a variety $X$ is split Tate as in \ref{THE THEOREM}, for example if $X$ is a smooth projective variety admitting an action of $\G_m$ (see \cite{Brosnan}). Similarly, \cite{Williams} shows that the motive of $GL_n(K)$ is split Tate. However, in the quasi-projective setting, one might ask for more general instances of such splittings, perhaps when the variety in question is not even a group scheme, such as $\Sym^{n,k}$.\\
    A first approach to computing the cohomology of $\Sym^{n,k}$ comes from understanding the inductive point-counting formula, which was elegantly computed for all $n,k\in\N$ (and for \emph{all} primes $p$) in \cite{M} via ``a heroic piece of elementary algebra" and in a way such that ``the most sophisticated mathematical notion used is that the row rank of a matrix equals the column rank." One can promote this formula to a statement about the isomorphism class of $[\Sym^{n,k}]$ in the \emph{Grothendieck ring of varieties}, $K_0(\text{Var}_K)$. This is the free abelian group generated by isomorphism classes $[X]$ of varieties over $K$ with relations $[X]=[Z]+[X - Z]$ for any closed embedding $Z\hookrightarrow X$, and with ring structure induced by the Cartesian product of varieties. Upon further observation, the point counting formula in \cite{M} suggests that we construct a filtration (see \ref{thefiltration}) on  $\Sym^{n,\leq k}$ over any field which is given by the rank of the (1,1) minor, and recover $[\Sym^{n,k}]$ via the graded pieces. In Section 2, we show the following:
\begin{theorem}
Fix any field $K$ of characteristic not $2$. In $K_0(\mathrm{Var}_K)$, and for $k>n$ or $k<0$, $[\Sym^{n,k}]=0$. Further, for all $n\in\N$, $[\Sym^{n,0}]=1$, and for $0< k\leq n$, $[\Sym^{n,k}]$ is the following polynomial in $\L$:\begin{footnote}{We are not inverting $\L$, but are writing the formulas in a more presentable manner.}
\end{footnote}
    \begin{equation}\label{main:thm}
    [\Sym^{n,k}]=\begin{cases}
        \Big(\prod_{i=1}^{k/2}\frac{\L^{2i}}{\L^{2i}-1}\Big)\Big(\prod_{i=0}^{k-1}(\L^{n-i}-1)\Big)\quad&\text{ if $k$ is even}\\
        \Big(\prod_{i=1}^{(k-1)/2}\frac{\L^{2i}}{\L^{2i}-1}\Big)\Big(\prod_{i=0}^{k-1}(\L^{n-i}-1)\Big)\quad&\text{ if $k$ is odd}
    \end{cases}
    \end{equation}
    where $\L:=[\A^1]$ denotes the Lefschetz motive.
\end{theorem}
The lemmas of this section show that the (1,1)-minor projection map:
\[
\pi_k:\Sym^{n,\leq k}\to \Sym^{n-1,\leq k}
\]
restricted to each $\Sym^{n-1,\ell}\subset \Sym^{n-1,\leq k}$, for $\ell\leq k$, is an algebraic vector bundle. 
In section 3, we use the filtration in \ref{thefiltration} to prove \ref{THE THEOREM}. While this filtration is the key construction in the paper, we define and prove some general instances in which motives can be split proper Tate in $C^{gm}$, for example, see \ref{prop:motivecomplement}. 

From our analysis, one might believe that inductive point counting formulas always lead to such splittings, but this would be false. If the coefficient ring is not a field, for example $\Z$, the space of $n\times n$ projective symmetric matrices of full rank, $\P\Sym^{n,n}$, may admit a similar point counting formula, but for $n=2$ and $k=\C$, the variety $\P\Sym^{2,2}$ cannot be split Tate as $H^2(\P\Sym^{2,2}(\C),\Z)\cong \Z/2\Z$. However, in our setting, it is in fact split Tate under certain conditions on the coefficient field (see \ref{projectivization}). There are further examples of singular proper Toric varieties, and more generally, \emph{linear schemes} in the sense of \cite{Totaro} which are not split Tate but are inductively defined and therefore allow for an inductive point counting formula over finite fields (see Section 7 of \cite{Totaro} for a concrete example). Regardless, such inductive point counting formulas shed significant light on the motives in question, either forcing them to split, or allowing for a more explicit computation of the extensions they define. \\\\

\textbf{Acknowledgements.} The author acknowledges Jesse Wolfson, his advisor, for countless insightful discussions and suggestions. The author thanks Burt Totaro whose remarks pointed out an error in the previous draft, significantly improving the paper. The author also thanks Francis Brown for suggesting the question. 
\section{Computations in the Grothendieck Ring}
Without doing much, we already have the following lemma:
\begin{lemma}
In $K_0(\text{Var}_K)$, 
\[[\Sym^{n,k}]=[\Sym^{n,\leq k}]-[\Sym^{n,\leq k-1}],\]
and in particular,
\[[\Sym^{n,n}]=\L^{n(n+1)/2}-[\Sym^{n,\leq n-1}]\ .\]
\end{lemma}
\begin{proof}
    First, observe that for all $k$, $\Sym^{n,\leq k}$ is a closed subvariety of $\Sym^{n,\leq n}\cong\A^{n(n+1)/2}$, as it is given by the intersection of the vanishing loci of the determinant polynomials of all $(k+1)\times (k+1)$ minors. Further, the complement of the closed embedding
    \[\Sym^{n,\leq k-1}\hookrightarrow \Sym^{n,\leq k}\]
    is precisely $\Sym^{n,k}$, giving us the desired cut and paste relations in $K_0(\text{Var}_K)$.
\end{proof}
We now define the filtration discussed in the introduction.

\begin{definition}
Define $\pi_k:\Sym^{n,\leq k}\to \Sym^{n-1,\leq k}$ as the morphism sending a matrix $(x_{ij})_{i,j=1,\dots,n}$ to its $(1,1)$ minor $(x_{ij})_{i,j=2,\dots,n}$.
\end{definition}
\begin{proposition}\label{thefiltration}
    For $k<n$, we can define an increasing filtration of length $k+1$:
    \[F_{\ell}(\Sym^{n,\leq k}):=\pi_k^{-1}(\Sym^{n-1,\leq \ell})\] which gives us a diagram in $\text{Var}_K$:
    \begin{equation}\label{filtration}
    \begin{tikzcd}
        \O\arrow[r,hook]&F_{0}(\Sym^{n,\leq k})\arrow[r,hook]&\dots\arrow[r,hook]&F_{k-1}(\Sym^{n,\leq k})\arrow[r,hook]& F_k(\Sym^{n,\leq k})=\Sym^{n,\leq k}\\
        &\O\arrow[r,hook]\arrow[u,"\circ",swap]&\dots\arrow[r,hook]&\pi_k^{-1}(\Sym^{n-1,[1,k-1]})\arrow[u,"\circ",swap]\arrow[r,hook]&\pi_k^{-1}(\Sym^{n-1,[1,k]})\arrow[u,"\circ",swap]\\
        &&&\vdots\arrow[u,"\circ",swap]&\vdots\arrow[u,"\circ",swap]\\
        &&&\pi_k^{-1}(\Sym^{n-1,k-1})\arrow[u,"\circ",swap]\arrow[r,hook]&\pi_k^{-1}(\Sym^{n-1,[k-1,k]})\arrow[u,"\circ",swap]\\
        &&&\O\arrow[r,hook]\arrow[u,"\circ",swap]&\pi_k^{-1}(\Sym^{n-1,k})\arrow[u,"\circ",swap]\\
        &&&&\O\arrow[u,"\circ",swap]
    \end{tikzcd}
    \end{equation}
    Here $\Sym^{n,[k,\ell]}$ denotes the variety of $n\times n$ symmetric matrices of rank $m$ such that $k\leq m\leq \ell$. The horizontal morphisms are closed embeddings, the vertical morphisms are open immersions, and if $X_{ij}$ denotes the variety in the $i$'th row and $j$'th column in \ref{filtration}, then the scheme-theoretic image of $X_{jk}\to X_{ik}$ is the complement of the closed embedding $X_{ij}\to X_{ik}$.
\end{proposition}
\begin{proof}
Since $\Sym^{n-1,\leq \ell}\hookrightarrow \Sym^{n-1,\leq \ell+1}$ is a closed embedding, we have that $\pi_k^{-1}(\Sym^{n-1,\leq \ell})\to \pi_k^{-1}(\Sym^{n,\leq \ell-1})$ must also be a closed embedding, i.e. the first row of \ref{filtration} is a sequence of closed embeddings. If $X_{ij}$ denotes the variety in the i'th row and j'th column of \ref{filtration}, then it suffices to show that for $i\neq 1$, $X_{ij}$ is the complement of the closed embedding $F_i(\Sym^{n,\leq k})\hookrightarrow F_j(\Sym^{n,\leq k})$, and the rest of the statements follow. Observe that:
\begin{align*}
F_j(\Sym^{n,\leq k})-F_i(\Sym^{n,\leq k})&=\pi_k^{-1}(\Sym^{n-1,\leq j})-\pi_k^{-1}(\Sym^{n-1,\leq i})\\
&=\pi_k^{-1}(\Sym^{n-1,\leq j}-\Sym^{n-1,\leq i})\\
&=\pi_k^{-1}(\Sym^{n-1,[i,j]})
\end{align*}
which is precisely $X_{ij}$.
\end{proof}
\begin{corollary}\label{cor:recursion}

    In $K_0(\text{Var}_K)$:
    \[[\Sym^{n,\leq k}]=\Big[\pi_k^{-1}(\Sym^{n-1,\leq k-2})\Big]+\Big[\pi_k^{-1}(\Sym^{n-1,k-1})\Big]+\Big[\pi_k^{-1}(\Sym^{n-1,k})\Big]\]
\end{corollary}
Luckily, the projection map restricted to each piece of the right-hand side of \label{cor:2.4} is a Zariski vector bundle, as we will show in the following few lemmas. 
\begin{lemma}\label{lem:n-1,k-2}
 For all $k\leq n$, we have a pullback square
 \[
 \begin{tikzcd}
     \Sym^{n-1,\leq k-2}\times\A^n\arrow[r,hook]\arrow[d,"pr_1"]&\Sym^{n,\leq k}\arrow[d,"\pi_k"]\\
     \Sym^{n-1,\leq k-2}\arrow[r,hook]&\Sym^{n-1,\leq k}
 \end{tikzcd}
 \]
\end{lemma}
\begin{proof}
    Recall that if the rank of the (1,1) minor of a matrix is less than $k-2$, then the rank of the matrix must be less than $k$. Given $(n-1)\times (n-1)$ matrix $Y=(y_{ij})$ such that $rk(Y)\leq k-2$, then a point in the fiber $M\in \pi^{-1}_k(Y)$ is $n\times n$ symmetric matrix with $Y$ as a (1,1)-minor, and so $rk(M)\leq k$. Specifying a matrix in the fiber involves choosing the first row, i.e. any vector in $\A^n$, and so we have the isomorphism:
    \begin{align*}
    \Phi:\pi_k^{-1}(\Sym^{n-1,\leq k-2})&\xrightarrow{\cong} \Sym^{n-1,\leq k-2}\times\A^n\\
    (x_{ij})_{i,j=1,\dots,n}&\mapsto \Big((x_{ij})_{i,j=2,\dots,n},x_{11},\dots,x_{1n}\Big)
    \end{align*}
    proving the claim.
\end{proof}
\begin{corollary}\label{cor:2.6}
    In $K_0(\text{Var}_K)$ for $k\leq n$:
    \[\Big[\pi_k^{-1}(\Sym^{n-1,\leq k-2})\Big]=\L^n[\Sym^{n-1,\leq k-2}]\]
\end{corollary}
\begin{lemma}\label{lem:n-1,k}
    For all $k<n$, we have another pullback square:
    \[
    \begin{tikzcd}
        \gamma_{n-1,k}^*\arrow[d]\arrow[r,hook]&\Sym^{n,\leq k}\arrow[d,"\pi_k"]\\
        \Sym^{n-1,k}\arrow[r,hook]&\Sym^{n-1,\leq k}
    \end{tikzcd}
    \]
    Where $\gamma_{n-1,k}^*$ is the pullback of the tautological bundle $\gamma_{n-1,k}$ over 
    $Gr(n-1,k)$, the Grassmannian of $k$ planes in $\A^{n-1}$, along the map \[\Sym^{n-1,k}\to Gr(n-1,k)\] which sends a matrix to the subspace of $\A^{n-1}$ spanned by its row vectors.
\end{lemma}

\begin{proof}
    Given a symmetric matrix $X=(x_{ij})_{i,j=1\dots,n}$ of rank $k$ whose $(1,1)$-minor $(x_{ij})_{i,j=2\dots,n}$ is also of rank $k$, the vector $(x_{12},\dots, x_{1n})\in\A^{n-1}$ must be contained in the span of the row vectors of the (1,1)-minor, i.e. the span of $\{(x_{22},\dots,x_{2n});\dots; (x_{n2},\dots,x_{nn})\}$. So, there exists $\lambda_1,\dots, \lambda_{n-1}\in K$ such that 
    \begin{equation}\label{eqn:lincomb}
    (x_{12},\dots, x_{1n})=\lambda_{1}(x_{22},\dots,x_{2n})+\dots+\lambda_{n-1}(x_{n2},\dots,x_{nn})
    \end{equation}
    The first row of $X$, $(x_{11},x_{12},\dots, x_{1n})$, must also be in the span of the bottom rows of $X$, otherwise the rank of $X$ will be greater than $k$. Since $x_{ij}=x_{ji}$, it follows that
    \begin{equation}\label{eqnlincomb}
    x_{11}=\lambda_2x_{12}+\dots+\lambda_{n-1}x_{1n}
    \end{equation}
    where the $\lambda_i$ are the same as in \ref{eqn:lincomb}. Thus, choosing 
    an $(n-1)\times (n-1)$ matrix $(x_{ij})_{i,j=2,\dots, n}$ of rank $k$, along with a vector 
    \[(x_{12},\dots, x_{1n})\in Span (\{(x_{22},\dots,x_{2n});\dots; (x_{n2},\dots,x_{nn})\})\]
    will uniquely determine a rank $k$ symmetric matrix whose (1,1) minor is also rank $k$. We use Pl\"{u}cker coordinates to define the total space of the tautological bundle $\gamma_{n-1,k}$, i.e.
    \[\gamma_{n-1,k}:=\{(\vec{v}_1\wedge\dots\wedge\vec{v}_k,\vec{x})\in Gr(n-1,k)\times\A^{n-1}: \vec{v}_1\wedge\dots\wedge\vec{v}_k\wedge\vec{x}=0\}
    \] and consider an open cover of $\Sym^{n-1,k}$ as follows: for $I=\{i_0<\dots<i_k\}$, let
    \[
    U_I:=\Bigg\{\begin{bmatrix}
        \vec{v_1}\\
        \vdots\\
        \vec{v_{n-1}}
    \end{bmatrix},\ rank\Big(\{v_{i_1},\dots, v_{i_k}\}\Big)=k\Bigg\}
    \]
    This defines an open cover $\{\pi^{-1}_k(U_I)\}_I$ on $\pi^{-1}_k(\Sym^{n-1,k})$. For each $I=\{i_0<\dots<i_k\}\subset \{1,\dots,n-1\}$ there is a map $\Phi_I:\pi^{-1}_k(U_I)\to \gamma_{n-1,k}$ given by
    \[
    \Phi_I:\begin{bmatrix}
        x_{11}&\vec{x}\\
        \vec{x}& \begin{bmatrix}
        \vec{v_1}\\
        \vdots\\
    \vec{v_{n-1}}
    \end{bmatrix}
    \end{bmatrix}\mapsto( \vec{v_{i_1}}\wedge\dots\wedge \vec{v_{i_k}},\vec{x})
    \]
    The maps $\Phi_I$ agree on intersections $\pi_k^{-1}(U_I\cap U_{I'})$ and so extend to a map $\Phi:\pi^{-1}_k(\Sym^{n-1,k})\to\gamma_{n-1,k}$. For $I=\{i_0<\dots<i_k\}$ we can further restrict $U_I$ to:
    \[
    U_{IJ}:=\{M\in U_I: det(M_{IJ})\neq 0\}
    \]
    where $M_{IJ}$ is the $k\times k$ minor of $M$ with rows 
    $I=\{i_0<\dots<i_k\}$ and columns $J=\{j_0<\dots<j_k\}$. We now check on affine opens that for each $I,J$: 
    \[
    \begin{tikzcd}
        \pi^{-1}_k(U_{IJ})\arrow[rr,"\Phi"]\arrow[d]&&\gamma_{n-1,k}\arrow[d]\\
       U_{IJ}\arrow[rr,"\Phi_{1}"]&&Gr(n-1,k)
    \end{tikzcd}
    \]
    is a pullback diagram, where $\Phi_{1}(M):=\vec{v_{i_1}}\wedge\dots\wedge \vec{v_{i_k}}$. We can define a morphism
    \[\Psi:U_{IJ}\times_{Gr(n-1,k)}\gamma_{n-1,k}\to \pi^{-1}_k(U_{IJ})\]
    by sending \[\Big(M,\vec{v_1}\wedge\dots\wedge \vec{v_k},\vec{x}\Big)\mapsto 
    \begin{bmatrix}
    M_{IJ}\vec{x}\cdot \vec{x}&\vec{x}\\
    \vec{x}&M
    \end{bmatrix}\]
    Note that $M_{IJ}\vec{x}\cdot \vec{x}=\lambda_1x_{12}+\dots+\lambda_nx_{1n}$ as in \ref{eqnlincomb}. One can verify that this map is a two-sided inverse to the map $\pi^{-1}_k(U_{IJ})\to U_{IJ}\times_{Gr(n-1,k)}\gamma_{n-1,k}$ given by the universal property of the pullback. This concludes the proof.
\end{proof}
\begin{corollary}\label{cor:2.8}
    In $K_0(\text{Var}_K)$, we have $[\pi_k^{-1}(\Sym^{n-1,k})]=\L^{k}[\Sym^{n-1,k}]$
\end{corollary}
\begin{lemma}\label{lem:n-1,k-1}
For all $k<n$, we have a pullback square:
\[
\begin{tikzcd}
    \A^1\oplus\gamma_{n-1,k-1}\arrow[r,hook]\arrow[d]&\Sym^{n,\leq k}\arrow[d,"\pi_k"]\\
    \Sym^{n-1,k-1}\arrow[r,hook]&\Sym^{n-1,\leq k}
\end{tikzcd}
\]
Where $\gamma_{n-1,k-1}$ is the pullback of the tautological bundle over $Gr(k-1,n-1)$ along the map \[\Sym^{n-1,k-1}\to Gr(k,n-1)\] sending a matrix to the subspace of $\A^{n-1}$ spanned by its row vectors. 
\end{lemma}
\begin{proof}
    The proof is similar to that of \ref{lem:n-1,k}. We first prove by contradiction that if $X=(x_{ij})$ is symmetric of rank $k$ whose (1,1) minor $M$ is rank $k-1$, then the vector $(x_{12},\dots x_{1n})$ must be in the span of the rows of $M$. Suppose this is not the case. Then, the first column of $X$ will be a linear combination of the other columns, else $rk(X)>k$. However, $X$ is a symmetric matrix, which implies that the first row is a linear combination of the other rows, forcing $(x_{12},\dots x_{1n})$ to also be a linear combination of the rows of $M$ which contradicts the claim. So, there exists $\lambda_2,\dots,\lambda_{n-1}\in K$ such that
    
    \[
    (x_{12},\dots, x_{1n})=\lambda_{2}(x_{22},\dots,x_{2n})+\dots+\lambda_{n-1}(x_{n2},\dots,x_{nn})
    \] If $x_{11}=\sum_{i=2}^n\lambda_ix_{1i}$, then $rk(X)=k-1$ just as in the proof of \ref{lem:n-1,k}; otherwise, the rank of $X$ is $k$. So, a matrix in $\pi_k^{-1}(Y)$, where $Y\in \Sym^{n-1,k-1}$ is determined a vector $(x_{12},\dots,x_{1n})$ in the span of the rows of $Y$, along with any choice of $x_{11}\in K$. Just as in \ref{lem:n-1,k}, we have a cover $U_{I}$ of $\Sym^{n-1, k-1}$ but now $I=\{i_1<\dots< i_{k-1}\}$. We can then define $\Phi_{I}:\pi_k^{-1}(U_{I})\to \A^1\times \gamma_{n-1,k-1}$ sending:
    \[(x_{ij})\mapsto \Big(x_{11};v_{i_1}\wedge\dots\wedge v_{i_{k-1}},(x_{12},\dots,x_{1n})\Big)\]
    Exactly as in \ref{lem:n-1,k}, this extends to a map on $\pi_k^{-1}(\Sym^{n-1,k-1})$. Similar to \ref{lem:n-1,k}, we can conclude that we have a pullback square:
    \[
    \begin{tikzcd}
    \pi_k^{-1}(\Sym^{n-1,k-1})\arrow[r,"\Phi"] \arrow[d,"\pi_k"]&\A^1\oplus \gamma_{n-1,k-1}\arrow[d]\\
    \Sym^{n-1,k-1}\arrow[r,"\Phi_1"] &Gr(n-1,k-1)
    \end{tikzcd}
    \]
\end{proof}
\begin{corollary}\label{cor:2.10}
    In $K_0(\text{Var}_K)$, we have $[\pi_k^{-1}(\Sym^{n-1,k-1})]=\L^{k}[\Sym^{n-1,k-1}]$
\end{corollary}
Now, we can recover a recursive formula for the class of symmetric matrices of rank $k$:
\begin{corollary}\label{cor:2.11}
    We have that in $K_0(\text{Var}_K)$ for $0<k\leq n$:
    \[[\Sym^{n,k}]=(\L^n-\L^{k-1})[\Sym^{n-1,k-2}]+(\L^{k}-\L^{k-1})[\Sym^{n-1,k-1}]+\L^{k}[\Sym^{n-1,k}]
    \]
\end{corollary} 
\begin{proof}
   If $k<n$, then substituting the results from \ref{cor:2.6},\ref{cor:2.8}, and \ref{cor:2.10} into the equation in \ref{cor:recursion}, we get the relation
    \[
    [\Sym^{n,\leq k}]=\L^{n}[\Sym^{n-1,\leq k-2}]+\L^k[\Sym^{n-1,k-1}]+\L^k[\Sym^{n-1,k}]
    \]
   So for $k<n$:
    \begin{align*}
        [\Sym^{n,k}]&=[\Sym^{n,\leq k}]-[\Sym^{n,\leq k-1}]\\
        &=\L^{n}([\Sym^{n-1,\leq k-2}]-[\Sym^{n-1,\leq k-3}])+\L^{k}([\Sym^{n-1,k-1}]+[\Sym^{n-1,k}])-\L^{k-1}([\Sym^{n-1,k-2}]+[\Sym^{n-1,k-1}])\\
        &=(\L^n-\L^{k-1})[\Sym^{n-1,k-2}]+(\L^{k}-\L^{k-1})[\Sym^{n-1,k-1}]+\L^{k}[\Sym^{n-1,k}]
    \end{align*}
    If $k=n$, then $\Sym^{n,\leq n}$ is the affine space that represents all $n\times n$ symmetric matrices, and so the projection $\pi_n:\Sym^{n,\leq n}\to \Sym^{n-1,\leq n}$ is a trivial $n$-bundle, so (trivially)
\begin{align*}
    [\Sym^{n,\leq n}]&=\L^n\Big([\Sym^{n-1,\leq n-3}]+[\Sym^{n-1, n-2}]+[\Sym^{n-1, n-1}]\Big)
\end{align*}
Now we see that
\begin{align*}
    [\Sym^{n,n}]&=[\Sym^{n,\leq n}]-[\Sym^{n,\leq n-1}]\\
    &=(\L^n-\L^{n-1})([\Sym^{n-1,n-2}])+(\L^n-\L^{n-1})[\Sym^{n-1,n-1}]
\end{align*}
Which is the correct formula, as $\Sym^{n-1,n}=\O$. This concludes the proof.
\end{proof}
In order to get a general closed formula for $[\Sym^{n,k}]$ as a polynomial in $\L$, we use the same procedure as in theorem 2 of \cite{M}:
\begin{theorem}
    For $k>n$ or $k<0$, $[\Sym^{n,k}]=\O$. Further, for all $n\in\N$, $[\Sym^{n,0}]=1$, and for $0< k\leq n$ in $K_0(\mathrm{Var}_{K})$, $[\Sym^{n,k}]$ is the following polynomial in $\L$:
    \begin{equation}\label{main:thm}
    [\Sym^{n,k}]=\begin{cases}
        \Big(\prod_{i=1}^{k/2}\frac{\L^{2i}}{\L^{2i}-1}\Big)\Big(\prod_{i=0}^{k-1}(\L^{n-i}-1)\Big)\quad&\text{ if $k$ is even}\\
        \Big(\prod_{i=1}^{(k-1)/2}\frac{\L^{2i}}{\L^{2i}-1}\Big)\Big(\prod_{i=0}^{k-1}(\L^{n-i}-1)\Big)\quad&\text{ if $k$ is odd}
    \end{cases}
    \end{equation}
\end{theorem}
\begin{proof}
    We just need to show \ref{main:thm}. Since it is well-known that $\L$ is a zero-divisor in $K_0(\mathrm{Var}_K)$ (see \cite{Borisov}) we must first restrict ourselves to computations to the ring generated by $\L$ in $K_0(\mathrm{Var}_K)$. This is a free sub-ring, $\Z[\L]$, of $K_0(\mathrm{Var}_k)$ in which $\L$ is not a zero divisor (otherwise, $\L^n=0$ for some $n>0$), and so we can preform left/right cancellation in $\Z[\L]$. We can quickly check that $[\Sym^{0,0}]=1, [\Sym^{1,0}]=1,[\Sym^{1,1}]=[K^\times]=\L-1$ agree with the formula, so we just need to show that the formula satisfies the recursion relation in \ref{cor:2.11}. If $k$ is even, then write $f(\L)=\prod_{i=1}^{k/2}\frac{\L^{2i}}{\L^{2i}-1}\in \Z[\L,\L^{-1}]$. Plugging the formula into \ref{cor:2.11} and doing some re-indexing, we have that in $\Z[\L,\L^{-1}]$:
    \begin{align*}
      f(\L)\prod_{i=0}^{k-1}(\L^{n-i}-1)&= (\L^{n}-\L^{k-1})\Big(\prod_{i=1}^{(k-2)/2}\frac{\L^{2i}}{\L^{2i}-1}\Big)\Big(\prod_{i=1}^{k-2}\L^{n-i}-1\Big)\\
      &+(\L^k-\L^{k-1})\Big(\prod_{i=1}^{(k-2)/2}\frac{\L^{2i}}{\L^{2i}-1}\Big)\Big(\prod_{i=1}^{k-1}\L^{n-i}-1\Big)\\
      &+\L^k f(\L)\prod_{i=1}^{k}(\L^{n-i}-1)\\
      f(\L)\Big(\prod_{i=1}^{k-1}(\L^{n-i}-1)\Big)[(\L^{n}-1)-\L^k(\L^{n-k}-1)]&=\Big(\prod_{i=1}^{(k-2)/2}\frac{\L^{2i}}{\L^{2i}-1}\Big)\Big(\prod_{i=1}^{k-2}\L^{n-i}-1\Big)[(\L^{n}-\L^{k-1})+(\L^k-\L^{k-1})(\L^{n-k+1}-1)]\\
      f(\L)\Big(\prod_{i=1}^{k-1}(\L^{n-i}-1)\Big) (\L^k-1)&=\Big(\prod_{i=1}^{(k-2)/2}\frac{\L^{2i}}{\L^{2i}-1}\Big)\Big(\prod_{i=1}^{k-2}\L^{n-i}-1\Big)[\L^{n+1}-\L^k]\\
      f(\L)(\L^k-1)&=\L^{k}\prod_{i=1}^{(k-2)/2}\frac{\L^{2i}}{\L^{2i}-1}\\
      f(\L)&=\prod_{i=1}^{k/2}\frac{\L^{2i}}{\L^{2i}-1}
    \end{align*}
    This string of equalities shows us that the recursion formula is satisfied in $\Z[\L,\L^{-1}]$, and a similar proof works when $k$ is odd. Since the formula is a polynomial, it also holds $\Z[\L]$. Lastly, since $\Z[\L]$ is a sub-ring of $K_0(\mathrm{Var}_K)$, the formula holds there as well. 
\end{proof}
The formulas are a bit simpler if we want to compute the class of non-singular symmetric matrices:
\begin{corollary}
    For $n\in\N$:
    \begin{align*}
        [\Sym^{n,n}]=\begin{cases}
            \prod_{i=1}^{n/2}(\L^{n+1}-\L^{2i})\quad&\text{ if $n$ is even}\\
            \prod_{i=0}^{(n-1)/2}(\L^{n}-\L^{2i})\quad&\text{ if $n$ is odd}\\
        \end{cases}
    \end{align*}
\end{corollary} 
\section{Computing the Motive}
Let $K$ be any perfect field of characteristic not equal to $2$, and let $\F$ be any field. If $K$ admits resolution of singularities, we work in $DM^{eff}_{gm,Nis}(K,\F)$, the triangulated category of Voevodsky's effective geometric motives as in \cite{Vo1} with coefficients in $\F$, and if $K$ is perfect of characteristic $p$, then we can use $DM^{eff}_{gm,ldh}(K,\F)$, where $p\in \F$ is invertible, as in \cite{K}. Since the arguments in each setting are the same, we can just write $C^{gm}$ as the category we are working in.
\begin{definition}
In $C^{gm}$, let $\F(0):=M_{gm}(*)$ denote the motive of a point. The Tate motive, denoted $\F(1)$, is the complex $\tilde{M}_{gm}(\P^1)[-2]$, where for any variety $X$, $\tilde{M}_{gm}(X)$ is the complex in degrees $0$ and $-1$:
\[
M_{gm}(X)\to \F(0)
\]
where the morphism is induced by the terminal morphism $X\to \mathrm{Spec}(K)$. Define, for $n\geq 1$:
\[
\F(n):=\F(1)^{\otimes n}
\]
Similarly, for $A\in C^{gm}$ define $A(n):=A\otimes \F(n)$. Tensoring with $\F(1)$ is formally inverted in $C^{gm}$, and so we define $\F(-1)$ as its inverse and $\F(-n):=\F(-1)^{\otimes n}$ for $n>0$.
\end{definition}
Most of the next definition is directly from definition 2.60 of \cite{VoEilenberg}.
\begin{definition}
    We define an object to be split Tate if it is a direct sum $\F(n)[r]$, where $n,r\in\Z$. An object is split \emph{proper} Tate if it is a direct sum of $\F(n)[r]$, where $n\in\Z_{\geq 0}$ \emph{and} $r\geq 2m$. An object of $C^{gm}$ is split \emph{positive} Tate if it is a direct sum of $\F(n)[r]$, for $n\in\Z_{\geq 0}, r\in\Z$. Let $\mathrm{\overline{SPT}}$ (resp. $\mathrm{\overline{SDT}}$) denote the category of split proper (resp. positive) Tate motives in $C^{gm}$. Lastly, let $\mathrm{\overline{SDT}}_{\leq n}$ denote subcategory of split positive Tate motives with summands of the form $\F(a)[b]$ where $0<a<n$.   
\end{definition}
We will rely often on the following foundational results in $C^{gm}$:
\begin{theorem}[4.3.1 of \cite{Vo1} and 5.3.16 of \cite{K}]\label{invertibility}
    For any $A,B\in C^{gm}$, the map \[\mathrm{Hom}_{C^{gm}}(A,B)\to \mathrm{Hom}_{C^{gm}}(A(1),B(1))\] is an isomorphism. 
\end{theorem}
\begin{proposition}[4.1.8 of \cite{Vo1} and 5.3.9 of \cite{K}]\label{Anistate}
    In $C^{gm}$:
    \[
    M^c_{gm}(\A^n)\cong \F(n)[2n]
    \]
\end{proposition}
\begin{proposition}[4.1.5 of \cite{Vo1} and 5.3.5 of \cite{K}]\label{LES}
    Let $X$ be an algebraic variety, and $Z\subset X$ a closed subvariety. Then, we have an distinguished triangle of compactly supported motives in $C^{gm}$:
    \[
    M_{gm}^c(Z)\to M_{gm}^c(X)\to M_{gm}^c(X-Z)\to M_{gm}^c(Z)[1]
    \]
\end{proposition}

\begin{proposition}\label{cor:homvanish}
    $\mathrm{Hom}_{C^{gm}}(\F(0),\F(a)[b])=0$ if :
    \begin{enumerate}    
        \item $a=0,\ b\neq 0$
        \item $a<0$
        \item $b>1$
    \end{enumerate}
    Further, a morphism from $\F(a)[b]$ to $\F(a')[b']$ is invertible iff $a=a'$ and $b=b'$, where by \ref{invertibility}, $\mathrm{Hom}_{C^{gm}}(\F(a)[b],\F(a)[b])=\mathrm{Hom}_{C^{gm}}(\F(0),\F(0))=\F$.
\end{proposition}
\begin{proof}
Let $k$ denote the underlying field of $C^{gm}$ for this proof. For the first statement, we know that $\F(0)$ is a complex concentrated in degree $0$, which immediately implies $\mathrm{Hom}_{C^{gm}}(\F(0),\F(0)[b])=0$ for $b\neq 0$. If $b=0$, then by work of Voevodsky \cite{Vo1} 4.2.9 and \cite{K} 5.3.14, we have
\[
\mathrm{Hom}_{C^{gm}}(\F(0),\F(a)[b])=CH^{a-b}(\mathrm{Spec}(k),2a-b)_\F
\]
 Setting $b=0$ have the hom-set equal to $CH^{0}(\mathrm{Spec}(k),0)_\F=\mathrm{Ch}^0(\mathrm{Spec}(k))_\F=\F$, where the latter is the usual Chow ring. For the next two statements we use the identification by Bloch in \cite{Bloch}:
\[
CH^{a-b}(\mathrm{Spec}(k),2a-b)=\mathrm{gr}^a_\gamma K_{2a-b}(k)_\F
\]
where right-hand side denotes successive quotients of the decreasing $\gamma$-filtration $F^\bullet_\gamma K(k)$ defined in \cite{Soule} for \emph{connective} algebraic $K$ theory using the $\lambda$ operations. The second statement holds because the filtration is non-negative. For the third statement we see that for $b>1$, by the discussion in 1.5 of \cite{Soule}:
\[F^{a+1}_\gamma K_{2a-b}(k)\supset K_{2a-b}(k)\]
this implies that $\mathrm{gr}^a_\gamma K_{2a-b}(k)=0$ as desired. Lastly, if a morphism $\mathrm{Hom}_{C^{gm}}(\F(a)[b],\F(a')[b'])$ is invertible, then $\mathrm{Hom}_{C^{gm}}(\F(a)[b],\F(a')[b'])$ must also be non-zero, which can only happen if $a=a'$ by our analysis. In this case we have that the hom-set is equal to $\mathrm{Hom}_{C^{gm}}(\F(0),\F(0)[b'-b])$, which is non-zero iff $b'=b$, as desired. 
\end{proof}
\begin{definition}
     A morphism $f:A\to B$ in $C^{gm}(R)$ \emph{induces a splitting} if $ker(f),im(f)$ exist and are direct summands of $A$ and $B$ respectively, i.e. we have the diagram:
    \[
    \begin{tikzcd}
        A\arrow[r,"f"]\arrow[dd,"\cong",swap]&B\arrow[dd,"\cong"]\\
        &\\
        ker(f)\oplus im(f)\arrow[r,"\begin{bmatrix}
            0\ id\\
            0\ 0
        \end{bmatrix}"]& im(f)\oplus B'
    \end{tikzcd}
    \]
    
\end{definition}
\begin{lemma}\label{cor:coneSPT}
    If $A$ and $B$ are objects of $\mathrm{\overline{SPT}}$ or $\mathrm{\overline{SDT}}_{\leq n}$, a morphism $f:A\to B$ in $C^{gm}$ \emph{induces a splitting} and if and only if the cone of $f$ is also split proper Tate, i.e. given a distinguished triangle:
    \[
    A\xrightarrow{f} B\to C\to A[1]
    \]
    $C$ is also split proper Tate.
\end{lemma}
\begin{proof}
Assume $A,B\in\mathrm{\overline{SPT}}$ as the proof is the same for $\mathrm{\overline{SDT}}_{\leq n}$. If $f$ induces a splitting, then by 2.73 of \cite{VoEilenberg}, both are closed under direct summands, and so $ker(f), im(f),B'$ are direct summands of $A,B$ respectively making them split proper Tate. The cone of $f$ will also be split proper Tate and isomorphic to $B'\oplus ker(f)[1]$. Now, assume we have an exact triangle $A\xrightarrow{f} B\to C\to A[1]$ of split proper Tate motives. Then, we have:
\[A=\bigoplus_{a,i\geq 2a}\F(a)[i]^{\oplus m_1(a,i)},\quad\text B=\bigoplus_{a,i\geq 2a}\F(a)[i]^{\oplus m_2(a,i)}\]
where $m_1(a,i),m_2(a,i)\in \Z_{\geq 0}$ for all $a,i$. So any morphism from $A$ to $B$ is comprised uniquely by morphisms from $\F(a)[i]$ to $\F(b)[j]$ in the summand decomposition of $A$ and $B$. However, in order for the cone to be split proper Tate, the morphisms from $\F(a)[i]$ to $\F(b)[j]$ must be invertible or $0$ (otherwise the cones of these individual morphisms will not be split proper Tate). So, using \ref{cor:homvanish}, we can construct by-hand the image and kernel of $f$ as split proper Tate summands of $A$ and $B$ respectively.
\end{proof}
\begin{corollary}\label{cor:splitting}
    Given a morphism of distinguished triangles in $\mathrm{\overline{SPT}}$ or $\mathrm{\overline{SDT}}_{\leq n}$:
    \[\begin{tikzcd}        X\arrow[r,"g_1"]\arrow[d,"f_1"]&Y\arrow[r,"g_2"]\arrow[d,"f_2"]&Z\arrow[r,"+1"]\arrow[d,"f_3"]&X[1]\arrow[d,"f_1{[}1{]}"] \\        X'\arrow[r,"g_1{'}"]&Y'\arrow[r,"g_2{'}"]&Z'\arrow[r,"+1"]&X'[1] \\
    \end{tikzcd}
    \]
    If $f_1,f_3$ induce splittings, then so does $f_2$. 
\end{corollary}

\begin{proof}
Assume that the diagram is in $\mathrm{\overline{SPT}}$ (the proof is the same in $\mathrm{\overline{SDT}}_{\leq n}$). First observe that if $A,B\subset X$ are direct summands of $X\in \mathrm{\overline{SPT}}$ defined by projectors $p,q$ respectively, then $A\cap B:=A\times_X B=\mathrm{im}(p\circ q)$ is also a well-defined summand of $X$ in $\mathrm{\overline{SPT}}$ corresponding to the projector $p\circ q$ (which is $=q\circ p$ as the homs are abelian), as the images of all projectors are contained in $\mathrm{\overline{SPT}}$. The morphism of the triangles in the corollary then becomes
    
    \[
    \begin{tikzcd}
    A\oplus B\oplus C\oplus D\arrow[r]\arrow[d]&A\oplus B\oplus E\oplus F\arrow[r]\arrow[d]& E\oplus F\oplus A[1]\oplus B[1]\arrow[d]\\
    A\oplus B'\oplus C\oplus D'\arrow[r]& A\oplus B'\oplus E\oplus F'\arrow[r]&E\oplus F'\oplus C[1]\oplus D'[1]
    \end{tikzcd}
    \] 
    where 
    $D=\Big(ker(f_1)\cap ker( g_1)\Big)$, $C\oplus D=ker(g_1)$, $B\oplus D=ker(f_1)$, $C\oplus D'= ker(g_1')$ and $F=ker(f_3)$. The rest of the summands become identified with the various images of the morphisms, which are all split proper Tate. Further, all the morphisms are diagonal morphisms that are either the identity if the summand is the same, or $0$ otherwise. This shows that we can write $f_2=(id_A,0,id_E,0)$, giving us the induced splitting.
\end{proof}
\begin{corollary}\label{cor:compositionsplits}
    Let $g:A\to B$ and $f:B\to C$ be morphisms in $\mathrm{\overline{SPT}}$ (resp. $\mathrm{\overline{SDT}}_{\leq n}$). Then $f\circ g$ induces a splitting if both $f$ and $g$ induce splittings.
\end{corollary}
\begin{proof}
    In the diagram in corollary \ref{cor:splitting}, set $g_1=id$ and $f_1=f, g_1'=g$, and $f_2=f\circ g$. This forces $Z=0$, and $f_3=0$. Lastly, $Z'$ is in $\mathrm{\overline{SPT}}$ (resp. $\mathrm{\overline{SDT}}_{\leq n}$) by \ref{cor:coneSPT} as $g$ induces a splitting. The conditions in \ref{cor:splitting} apply, and we can conclude the result.
\end{proof}
One can already see the following standard result:
\begin{proposition}\label{n=1}
 $M_{gm}^c(\Sym^{1,1})$ is split proper Tate.
\end{proposition}
\begin{proof}
    Observe that $\Sym^{1,1}\cong\A^1-\{0\}$. The inclusion of the origin in $\A^1$ gives us an exact triangle:
    \[
    M^c_{gm}(*)\to M^c_{gm}(\A^1)\to M^c_{gm}(\A^{1}-\{0\})\to M^c_{gm}(*)[1]
    \]
    Using \ref{Anistate} and shifting gives us the exact triangle
    \begin{equation}\label{first0}
    \F(1)[2]\to M^c_{gm}(\A^{1}-\{0\})\to \F(0)[1]\to \F(1)[3]
    \end{equation}
    Now, we have
    \begin{align*}
        Hom_{C^{gm}}(\F(0)[1],\F(1)[3])&\cong Hom_{C^{gm}}(\F(0),\F(1)[2])\\
        &=0
    \end{align*}
    where the vanishing is a consequence of \ref{cor:homvanish}. So, \ref{first0} has a connecting homomorphism of $0$, giving us:
    \[
    M^c(\A^{1}-\{0\})\cong \F(0)[1]\oplus \F(1)[2]
    \]
    and so $M^c(\A^{1}-\{0\})$ is split proper Tate.
\end{proof}
Before the proof of the main theorem, we will need another proposition that is a quick consequence of theorems in \cite{K} and \cite{Vo1}:
\begin{proposition}\label{prop:motivetotalspace}
    If $\pi:E\to B$ is an algebraic vector bundle of rank $n$, then the induced map:
    \[
    \pi^*:M_{gm}^c(B)(n)[2n]\to M_{gm}^c(E)
    \]
    in 4.2.4 of \cite{Vo1} and 5.3.12 of \cite{K} is an isomorphism.
\end{proposition}
\begin{proof}
    We prove this by induction on the size of the trivializing chart $\{U_i\}_{\{i=1,\dots,k\}}$ of $\pi$. For $k=1$, this is the content of 4.1.8 of \cite{Vo1} and 5.3.9 of \cite{K}. Assuming the proposition is true for some $k$, we assume $\pi$ has a trivializing chart $\{U_k\}_{\{i=1,\dots,k+1\}}$. Writing $Z=X-U_{k+1}$, we see that $\{U_i\cap Z\}_{\{i=1,\dots,k\}}$ is an open cover for the bundle $\pi$ restricted to $Z$, i.e. $\pi^{-1}(Z)\to Z$. we have a diagram of distinguished triangles due to \ref{LES}:
    \[\begin{tikzcd}
        M_{gm}^c(\pi^{-1}(Z))\arrow[r]& M_{gm}^c(E)\arrow[r]& M_{gm}^c(\pi^{-1}(U_{k+1}))\arrow[r,"+1"]& \  \\
        M_{gm}^c(Z)(n)[2n]\arrow[r]\arrow[u,"\pi^*"]& M_{gm}^c(B)(n)[2n]\arrow[r]\arrow[u,"\pi^*"]& M_{gm}^c(U_{k+1})(n)[2n]\arrow[r,"+1"]\arrow[u,"\pi^*"]& \  
    \end{tikzcd}\]
    Since $\pi^{-1}(U)\cong U\times\A^n$, the right vertical map is an isomorphism, and the first vertical map is an isomorphism due to the inductive hypothesis. This forces $M_{gm}^c(B)(n)[2n]\to M_{gm}^c(E)$ to be an isomorphism as well, proving the lemma. 
\end{proof}
The next proposition is crucial to the main proof. 
\begin{proposition}\label{prop:motivecomplement}
    Let $E\to X$ be an $n$-dimensional vector bundle over a smooth equidimensional variety $X$, and let $F\to E$ be an $m$-dimensional sub-bundle. Then, in $C^{gm}$, if $M_{gm}^c(X)$ is split proper Tate, then so is the motive $M^c_{gm}(E\setminus F)$.
\end{proposition}
\begin{proof}
    First note that if $M^c_{gm}(X)$ is split proper Tate, we have a distinguished triangle:
    \[
    M^c_{gm}(F)\to M^c_{gm}(E)\to M^c_{gm}(E\setminus F)\xrightarrow{+1}M^c_{gm}(F)[1]
    \]
    where by \ref{prop:motivetotalspace}, $M^c_{gm}(F)\cong M^c_{gm}(X)(m)[2m]$ and  $M^c_{gm}(E)\cong M^c_{gm}(X)(n)[2n]$, making them split proper Tate. If we show that $M_{gm}(E\setminus F)$ is split Tate (i.e. in $\mathrm{\overline{SDT}}$), then it will follow by Poincar\'{e} duality in $C^{gm}$ (since $E\setminus F$ is smooth) that $M_{gm}^c(E\setminus F)$ is split Tate. Further, by \ref{cor:splitting}, $M^c_{gm}(E\setminus F)$ is the cone of a morphism of split proper Tate motives that induces a splitting, making in fact $M^c_{gm}(E\setminus F)$ split proper Tate. 
    
    So, it suffices to show that $M_{gm}(E\setminus F)$ is split Tate. We have the decomposition
    \[
    M^c_{gm}(X)= \bigoplus_{a\geq 0}\bigoplus_{i\geq 2a}\F(a)[i]^{\bigoplus m_i(a)}
    \]
    where $m_i(a)\geq 0$, and the sum is finite as $X$ is equidimensional of dimension $D$. By Poincar\'{e} duality (4.3.7 of \cite{Vo1} and 5.3.18 of \cite{K}):
    \[
    M_{gm}(X)=\bigoplus_{a\geq 0}\bigoplus_{i\leq 2a}\F(D-a)[2D-i]^{\bigoplus m_i(a)}
    \]
    Since the direct sum is finite as $X$ is a scheme of finite type, we can twist by some $\F(\ell)$ for $\ell\geq 0$ such that $M_{gm}(X)(\ell)\in \mathrm{\overline{SDT}}_{\leq D+\ell}$, and so without loss of generality, we assume  $M_{gm}(X)\in \mathrm{\overline{SDT}}_{\leq D}$, as the arguments will not change. By 4.18 of \cite{Deglise}, for an $n$-bundle $E\to X$ with an $m$-dimensional sub-bundle $F\to X$, we have a morphism of Gysin triangles:
    \[
    \begin{tikzcd}
    M_{gm}(F-X)\arrow[r]\arrow[d]& M_{gm}(X)\arrow[r]\arrow[d]& M_{gm}(X)(m)[2m]\arrow[r,"+1"]\arrow[d,"\Phi"]&M_{gm}(E-X)[1]\arrow[d]\\
    M_{gm}(E-X)\arrow[r]& M_{gm}(X)\arrow[r]& M_{gm}(X)(n)[2n]\arrow[r,"+1"]&M_{gm}(F-X)[1]
    \end{tikzcd}
    \]
    where $\Phi=\Big((1_{X})_*\boxtimes c_{n-m}(E/F)\Big)(m)[2m]$, which by definition is the composition
    \[
    \begin{tikzcd}
    M_{gm}(X)(m)[2m]\arrow[r,"M_{gm}(\Delta)(m){[}2m{]}"] &\Big(M_{gm}(X)\otimes M_{gm}(X)\Big)(m)[2m]\arrow[d,"(1_X)_*(m){[}2m{]}\otimes c_{n-m}(E/F)(m){[}2m{]}"]\\
    &\Big(M_{gm}(X)(m)[2m]\otimes\F(n-m)[2n-2m]\Big)(m)[2m]
    \end{tikzcd}
    \]
    Here $\Delta:X\to X\times X$ denotes the inclusion onto the diagonal, and $c_{n-m}(E/F):M_{gm}(X)\to \F(n-m)[2n-2m]$ is the Chern class map associated to the quotient bundle $E/F$. Now, the projection onto either coordinate $X\times X\to X$ is a right inverse, and so $M^{gm}(\Delta)$  (and therefore $M^{gm}(\Delta)(m)[2m]$) is right invertible, implying that it induces a splitting. Note that if $c_{n-m}(E/F)$ induces a splitting, then since $(1_{X})_*$ does so as well, $(1_{X})_*\otimes c_{n-m}(E/F)$ must also induce a splitting, and similarly  $\Big((1_{X})_*\boxtimes c_{n-m}(E/F)\Big)(m)[2m]$. 
    
    We have thus reduced the problem to showing that $c_{n-m}(E/F)$ induces a splitting. Set $n-m=d$, and observe that any morphism in
    \[
    Hom_{C^{gm}}(M_{gm}(X),\F(d)[2d])
    \]
    is uniquely determined by elements in $Hom_{C^{gm}}(\F(D-a)[2D-i],\F(d)[2d])$, where $\F(D-a)[2D-i]$ appears in the direct sum decomposition of $M_{gm}(X)$. Now, by \ref{invertibility},
    \begin{align*}
    Hom_{C^{gm}}(\F(D-a)[2D-i],\F(d)[2d])=Hom_{C^{gm}}(\F(0),\F(d+a-D)[2d+i-2D])
\end{align*}
By \ref{cor:homvanish} the right-hand side will be zero if $d+a-D<0$, we have that the only non-zero case can occur when both $d+a-D\geq 0$ and therefore $2d+i-2D\geq 0$ as $i\geq 2a$. Next, by \ref{cor:homvanish}, non-zero cases occur when $2d+i-2D=0,1$. If $2d+i-2D=1$, then it also follows by \ref{cor:homvanish} that $d+a-D\geq 1$, i.e. $2d+2a-2D\geq 2$. Since $i\geq 2a$, we have $2d+i-2D\geq 2$, a contradiction. So, the only possible case is when $2d+i-2D=0$, so $i=2D-2d$. In this case $D-a=d$ as well. So, the Hom set becomes
\[
Hom_{C^{gm}}(\F(d)[2d],\F(d)[2d])=\F
\]
This analysis shows that the only non-zero morphisms from any summand of $M_{gm}(X)\to \F(c)[2c]$ must come from summands of the form $\F(d)[2d]\subset M_{gm}(X)$. Now, let $\mathrm{\overline{SDT}}_{d,2d}$ denote the subcategory $\mathrm{\overline{SPT}}_{\leq D}$ with summands only of the form $\F(d)[2d]$. This is well-defined, as by 2.73 of \cite{VoEilenberg}, $\mathrm{\overline{SDT}}_{\leq D}$ contains all direct summands. We see that the functor $Hom_{C^{gm}}(-,\F(d)[2d])$ gives an isomorphism from $\mathrm{\overline{SDT}}_{d,2d}$ to the category of $\F-$vector spaces. The remark in 2.63 of \cite{VoEilenberg} applies in our context, and we have a cohomological functor $s_{d,2d}:\mathrm{\overline{SDT}}_{\leq D}\to \mathrm{\overline{SDT}}_{d,2d}$ sending a split positive Tate motive to its direct summands of the form $\F(d)[2d]$, and in our case the map $c_d(E/F): M_{gm}(X)\to \F(d)[2d]$ gives the following decomposition in $\mathrm{\overline{SDT}}$:
\[\begin{tikzcd}
s_{d,2d}(M_{gm}(X))\arrow[r,hook]\arrow[d,"\Phi"]&M_{gm}(X)\arrow[d,"c_d(E)"]\arrow[l,two heads, bend right =15,"\pi",swap]\\
   \F(d)[2d]\arrow[r,equal]&\F(d)[2d]
\end{tikzcd}
\]
Where $\pi$ denotes projection onto the summand, and $\Phi$ induces a splitting as it is in a morphism in $\mathrm{\overline{SDT}}_{d,2d}$ which is semi-simple, so $\mathrm{ker}(\Phi)\oplus \mathrm{im}(\Phi)\cong s_{d,2d}(M_{gm}(X))$  So, \[\mathrm{ker}(c_{d}(E/F))=\mathrm{ker}(\pi)\oplus \mathrm{ker}(\Phi)\] and
\[\mathrm{im}(c_{d}(E/F))=im(\pi)\]
These are both direct summands of $M_{gm}(X)$, implying that $c_{d}(E/F)$ induces a splitting, and that $M_{gm}(E-F)$ is split Tate. 
\end{proof}

\begin{remark}
        Note that \ref{prop:motivecomplement} fails when the coefficient ring is not a field, because $\mathrm{\overline{SDT}}_{d,2d}$ is not necessarily semi-simple anymore.
\end{remark}
\begin{lemma}\label{lem:[k-1,k]}
    Recall that $\Sym^{n,[k-1,k]}$ denotes the space of symmetric matrices of rank both $k-1$ and $k$. Then, the projection map: $\pi_{k}:\pi^{-1}_k(\Sym^{n,[k-1,k]})\to \Sym^{n,[k-1,k]}$ induces an isomorphism:
    \[
    M^c_{gm}(\pi^{-1}_k(\Sym^{n,[k-1,k]}))\cong M^c_{gm}(\Sym^{n,[k-1,k]})(k)[2k]
    \]
\end{lemma}
\begin{proof}
    Combining \ref{lem:n-1,k} and \ref{lem:n-1,k-1} with \ref{prop:motivetotalspace}, the projection maps induce isomorphisms
    \[
    M^c_{gm}(\pi^{-1}_k(\Sym^{n,k-1}))\xrightarrow{\cong} M^c_{gm}(\Sym^{n,k-1})(k)[2k]
    \]
    and 
    \[
    M^c_{gm}(\pi^{-1}_k(\Sym^{n,k}))\xrightarrow{\cong} M^c_{gm}(\Sym^{n,k})(k)[2k]
    \]
    Next, the open-closed decompositions \[\Sym^{n,k-1}\hookrightarrow \Sym^{n,[k-1,k]}\leftarrow \Sym^{n,k}\] and \[\pi_k^{-1}(\Sym^{n,k-1})\hookrightarrow \pi_k^{-1}(\Sym^{n,[k-1,k]})\leftarrow \pi_k^{-1}(\Sym^{n,k})\] give us a commuting diagram in $C^{gm}$:
    \[\begin{tikzcd}
        M^c_{gm}(\pi^{-1}(\Sym^{n,k-1}))\arrow[r]& M^c_{gm}(\pi_k^{-1}(\Sym^{n,[k-1,k]}))\arrow[r]&M^c_{gm}(\pi^{-1}(\Sym^{n,k-1}))\arrow[r]&M^c(\pi^{-1}(\Sym^{n,k-1}))[1]\\
        M^c(\Sym^{n,k-1})(k)[2k]\arrow[r]\arrow[u,"\pi^*_k"]&M^c(\Sym^{n,[k-1,k]})(k)[2k]\arrow[r]&M^c(\Sym^{n,k})(k)[2k]\arrow[u,"\pi^*_k"]\arrow[r]&M^c(\Sym^{n,k-1})(k)[2k+1]\arrow[u,"\pi^*_k{[}1{]}"]
    \end{tikzcd}\]
    Using \ref{lem:n-1,k} and $\ref{lem:n-1,k-1}$, we can invoke \ref{prop:motivetotalspace} to say that each map $\pi^*$ is an isomorphism. Further, just by elementary properties of triangulated categories, there exists an isomorphism
    \[
    M^c_{gm}(\Sym^{n,[k-1,k]})(k)[2k]\xrightarrow{\cong} M^c_{gm}(\pi_k^{-1}(\Sym^{n,[k-1,k]}))
    \]
    which completes the diagram, giving us the result.
\end{proof}
\begin{theorem}
    For all $n,k$, $M^c_{gm}(\Sym^{n,\leq k})$ is split proper Tate, and the map:
    \[
    M^c_{gm}(\Sym^{n,\leq k-1})\to M^c_{gm}(\Sym^{n,\leq k})
    \]
    induces a splitting, implying by \ref{cor:compositionsplits} that for 
    $k\leq \ell$:
    \[
    M^c_{gm}(\Sym^{n,\leq k})\to M^c_{gm}(\Sym^{n,\leq \ell})
    \]
    also induce splittings. Thus, the compactly supported motives $M^c_{gm}(\Sym^{n,k}),M^c_{gm}(\Sym^{n,[k,\ell]})$ are split proper Tate in $C^{gm}$.
\end{theorem}
\begin{proof}
    We prove this by induction. For $n=1$, we have already proved the only non-trivial case in \ref{n=1}. If $k=0$ or $n$ then note that $\Sym^{n,\leq k}$ is affine space, so it is immediately split proper Tate. Next, $\Sym^{n,\leq k}=\O$ for $k<0$ or $k>n$. If $0<k<n$, our goal is to use the diagram of varieties:
    \[
    \begin{tikzcd}
        \pi^{-1}_{k-1}(\mathrm{Sym}^{n-1,\leq k-2})\arrow[r,hook]\arrow[d,hook]&\mathrm{Sym}^{n,\leq k-1}\arrow[d,hook]&\pi^{-1}_{k-1}(\mathrm{Sym}^{n-1,k-1})\arrow[l,"\circ"]\arrow[d,hook]\\
        \pi_{k}^{-1}(\mathrm{Sym}^{n-1,\leq k-2})\arrow[r,hook]&\mathrm{Sym}^{n\leq k}&\pi_k^{-1}(\mathrm{Sym}^{n-1,[k-1,k]})\arrow[l,"\circ"]
    \end{tikzcd}
    \]
    Which gives us a morphism of distinguished triangles in $C^{gm}$ by \ref{LES}:
    \begin{equation}\label{masterdiagram}
    \begin{tikzcd}
        M^c_{gm}(\pi^{-1}_{k-1}(\mathrm{Sym}^{n-1,\leq k-2}))\arrow[r]\arrow[d,"f_1"]&M^c_{gm}(\mathrm{Sym}^{n,\leq k-1})\arrow[d,"f_2"]\arrow[r]&M^c_{gm}(\pi^{-1}_{k-1}(\mathrm{Sym}^{n-1,k-1}))\arrow[d,"f_3"]\arrow[r,"+1"]& \ \\
        M^c_{gm}(\pi_{k}^{-1}(\mathrm{Sym}^{n-1,\leq k-2}))\arrow[r]&M^c_{gm}(\mathrm{Sym}^{n\leq k})\arrow[r]& M^c_{gm}(\pi_k^{-1}(\mathrm{Sym}^{n-1,[k-1,k]}))\arrow[r,"+1"]& \ 
    \end{tikzcd}
    \end{equation}
    Recall that $\pi_{k}^{-1}(\mathrm{Sym}^{n-1,\leq k-2})\cong \mathrm{Sym}^{n-1,\leq k-2}\times\A^n$ by \ref{lem:n-1,k-2}. Further, we can use \ref{lem:n-1,k-1}, and \ref{lem:[k-1,k]} to see that the morphism in in \ref{masterdiagram} becomes
    \[
    \begin{tikzcd}
        M^c_{gm}(\pi^{-1}_{k-1}(\mathrm{Sym}^{n-1,\leq k-2}))\arrow[r]\arrow[d,"f_1"]&M^c_{gm}(\mathrm{Sym}^{n,\leq k-1})\arrow[d,"f_2"]\arrow[r]&\textcolor{red}{M^c_{gm}(\mathrm{Sym}^{n-1,k-1})(k-1)[2k-2]}\arrow[d,"f_3"]\arrow[r,"+1"]& \ \\
        \textcolor{red}{M^c_{gm}(\mathrm{Sym}^{n-1,\leq k-2})(n)[2n]}\arrow[r]&M^c_{gm}(\mathrm{Sym}^{n\leq k})\arrow[r]& \textcolor{red}{M^c_{gm}(\mathrm{Sym}^{n-1,[k-1,k]})(k)[2k]}\arrow[r,"+1"]& \ 
    \end{tikzcd}
    \]
    By the inductive hypotheses, the motives red in color are split proper Tate. The rest of the proof goes as follows:
    \begin{enumerate}
        \item We prove $M^c_{gm}(\pi_{k-1}^{-1}(\Sym^{n-1,\leq k-2}))$ is split proper Tate for all $0<k<n$, and that the morphism $f_1$ induces a splitting. 
        \item We prove that the connecting homomorphism of the bottom triangle of \ref{masterdiagram} induces a splitting for all $0<k<n$, implying that $M^c_{gm}(\Sym^{n,\leq k})$ is split proper Tate for all $0<k<n$. 
        \item We prove that $f_3$ induces a splitting, and invoke \ref{cor:splitting} to conclude.
    \end{enumerate}
    \begin{enumerate}
    \item[Proof of (1):] We have the following diagram of varieties:
    \begin{equation}\label{k-2k-1}
    \begin{tikzcd}
    \pi_{k-1}^{-1}(\Sym^{n-1,\leq k-3})\arrow[r,hook]\arrow[d,equal]& \pi_{k-1}^{-1}(\Sym^{n-1,\leq k-2})\arrow[d,hook]&\pi_{k-1}^{-1}(\Sym^{n-1,k-2})\arrow[l,"\circ"]\arrow[d,hook]\\
    \Sym^{n-1,\leq k-3}\times\A^n\arrow[r,hook]&\Sym^{n-1,\leq k-2}\times\A^n&\Sym^{n-1,k-2}\times\A^n\arrow[l,"\circ"]
    \end{tikzcd}
    \end{equation}
    where the vertical arrows denote the inclusions into the various spaces of symmetric matrices with no restrictions on the top row. This in turn gives us a morphism of exact triangles in $C^{gm}$. The first vertical arrow from the left gives equality by \ref{lem:n-1,k-2}, and $\pi_{k-1}^{-1}(\Sym^{n-1,k-2})$ is the total space of a vector bundle of rank $k-1$ over $\Sym^{n-1,k-2}$. When looking at the induced morphisms of distinguished triangles given by the rows of \ref{k-2k-1}, we have that each morphism in the bottom row must induce a splitting by \ref{cor:coneSPT} (because all the motives in the bottom row are split proper Tate by the inductive hypothesis). Further, we have a map of connecting homomorphisms:
    \[\begin{tikzcd}
        M^c_{gm}(\Sym^{n-1,k-2})(k-1)[2(k-1)]\arrow[r]\arrow[d,"g"]& M^c_{gm}(\Sym^{n-1,\leq k-3})(n)[2n+1]\\
        M^c_{gm}(\Sym^{n-1,k-2})(n)[2n]\arrow[ur,"h"]&
    \end{tikzcd}\]
    By the inductive hypothesis and by \ref{cor:coneSPT}, $h$ induces a splitting, and by \ref{prop:motivecomplement}, $g$ induces a splitting because it is the inclusion of the total space of a vector sub-bundle into larger (trivial) bundle so we can invoke \ref{prop:motivecomplement}. By \ref{cor:compositionsplits}, have that the composition $g\circ h$, i.e. the connecting homomorphism in the top row, must also induce a splitting implying that $M^c_{gm}(\pi_{k-1}^{-1}(\Sym^{n-1,\leq k-2}))$ is split proper Tate. To conclude, we now have a diagram in $C^{gm}$:
\[
    \begin{tikzcd}
    M^c_{gm}(\Sym^{n-1,\leq k-3})(n)[2n]\arrow[r,hook]\arrow[d,equal]& M^c_{gm}(\pi_{k-1}^{-1}(\Sym^{n-1,\leq k-2}))\arrow[d,"f_1"]\arrow[r]&M^c_{gm}(\Sym^{n-1,k-2})(k-1)[2(k-1)]\arrow[r,"+1"]\arrow[d,"g"]&\ \\
    M^c_{gm}(\Sym^{n-1,\leq k-3})(n)[2n]\arrow[r,hook]&M^c_{gm}(\Sym^{n-1,\leq k-2})(n)[2n]\arrow[r]&M^c_{gm}(\Sym^{n-1,k-2})(n)[2n]\arrow[r,"+1"]&\ 
    \end{tikzcd}    
\]
    where all the motives are split proper Tate, and the outer vertical morphisms induce splittings. By \ref{cor:splitting}, $f_1$ also must induce a splitting.
    
    \item[Proof of (2):] Similarl to part 1 have the following commutative diagram: \[
    \begin{tikzcd}
        \pi^{-1}_k(\Sym^{n-1,\leq k-2})\arrow[d,equal]\arrow[r,hook]&\Sym^{n,\leq k}\arrow[d,hook]&\pi^{-1}_k(\Sym^{n-1,[k-1,k]})\arrow[l,"
        \circ"]\arrow[d,hook]\\
        \Sym^{n-1,\leq k-2}\times\A^n\arrow[r]&\Sym^{n-1,\leq k}\times\A^n&\Sym^{n-1,[k-1,k]}\times\A^n \arrow[l]
    \end{tikzcd}
    \]
    where the vertical morphisms denote inclusions into the space of matrices of the respective minors with no restrictions on the top row. Similar to last time, using \ref{lem:[k-1,k]} we have a commutative diagram of connecting homomorphisms:
\[\begin{tikzcd}
        M^c_{gm}(\Sym^{n-1,[k-1,k]})(k)[2k]\arrow[r]\arrow[d,"g"]& M^c_{gm}(\Sym^{n-1,\leq k-2})(n)[2n+1]\\
        M^c_{gm}(\Sym^{n-1,[k-1,k]})(n)[2n]\arrow[ur,"h"]&
    \end{tikzcd}\]
    By the inductive hypothesis and by \ref{cor:coneSPT}, $h$ induces a splitting. In order to see that $g$ induces a splitting, we have $\pi^{-1}_k(\Sym^{n-1,k-1})$ is a subbundle of the trivial bundle $\Sym^{n-1,k-1}\times\A^n$ and similarly $\pi^{-1}_k(\Sym^{n-1,k})$ is a subbundle of $\Sym^{n-1,k}\times\A^n$, and the compatible open-closed decompositions of $\Sym^{n-1,[k-1,k]}\times\A^n$ and $\pi_k^{-1}(\Sym^{n-1,[k-1,k]})$ used in \ref{lem:[k-1,k]} give us
    \[
    \begin{tikzcd}
    M^c_{gm}(\Sym^{n-1,k-1})(k)[2k]\arrow[r]\arrow[d,"f_1"]& M^c_{gm}(\Sym^{n-1, [k-1,k]})(k)[2k]\arrow[d,"g"]\arrow[r]&M^c_{gm}(\Sym^{n-1,k})(k)[2k]\arrow[d,"f_3"]\\
    M^c_{gm}(\Sym^{n-1,k-1})(n)[2n]\arrow[r]&M^c_{gm}(\Sym^{n-1,[k-1,k]})(n)[2n]\arrow[r]&M^c_{gm}(\Sym^{n-1,k})(n)[2n]
\end{tikzcd}
\]
    By \ref{prop:motivecomplement}, $f_1,f_3$ induce splittings, and by \ref{cor:splitting}, so must $g$. Since both $g$ and $h$ induce splittings, so must $g\circ h$ by \ref{cor:compositionsplits}. This proves by \ref{cor:coneSPT} that the cone $g\circ h$ is split proper Tate, i.e. $M^c_{gm}(\Sym^{n,\leq k})[1]$, and therefore $M^c_{gm}(\Sym^{n,\leq k})$, for all $0<k<n$.

    \item[Proof of (3):] The morphism $f_3$ is the composition:
    \[
    M^c_{gm}(\pi_{k-1}^{-1}(\Sym^{n-1,k-1}))\to M^c_{gm}(\pi_{k}^{-1}(\Sym^{n-1,k-1}))\to M^c_{gm}(\pi_{k}^{-1}(\Sym^{n-1,[k-1,k]}))
    \]
The first morphism is induced by an inclusion of a sub-bundle into a vector bundle. This induces a splitting by \ref{prop:motivecomplement}, and the cone of the second morphism is $M^c_{gm}(\pi_{k}^{-1}(\Sym^{n-1,k}))$ which is split proper Tate by the inductive hypothesis and \ref{prop:motivetotalspace}, implying that it also must induce a splitting by \ref{cor:splitting}. Therefore, $f_3$ is the composition of morphisms that both induce splittings, so it must must also do so by \ref{cor:compositionsplits}. 
    \end{enumerate}
\end{proof}
Lastly, we have the following corollaries:
\begin{corollary}
    By Poincar\'{e} duality, $M_{gm}(\Sym^{n,k})$ is split Tate for all $n,k$. 
\end{corollary}

\begin{corollary}\label{inductive}
    The inclusion $\Sym^{n-1,n-1}\hookrightarrow\Sym^{n,n}$ given by sending:
    \[
    M\mapsto 
    \begin{bmatrix}
        1&0\\
        0&M
    \end{bmatrix}
    \]
    induces the inclusion of $M_{gm}(\Sym^{n-1,n-1})\otimes M_{gm}(\G_m)$ as a direct summand of $\Sym^{n,n}$
\end{corollary}
\begin{proof}
    For this proof, let $\pi_{n,n}^{-1}(\Sym^{n-1,k})$ denote the subspace of symmetric $n\times n$ matrices of full rank whose $(1,1)$ minor is of rank $k$. We first show that
    \[
    \Sym^{n,n}=\pi_{n,n}^{-1}(\Sym^{n-1,n-1})\coprod \pi_{n,n}^{-1}(\Sym^{n-1,n-2})
    \]
    This is because $\pi_{n,n}^{-1}(\Sym^{n-1,n-1})$ is both an open and closed subset of $\Sym^{n,n}$. It is open, as it is the complement of the determinant polynomial of the $(1,1)$ minor. For $M\in \Sym^{n,n}$ we write 
    \[
    M=\begin{bmatrix}
        \alpha&\vec{v}\\
        \vec{v}&M_{(1,1)}
    \end{bmatrix}
    \]
    If the rank $M_{(1,1)}$ is full, $\vec{v}$ must be in the span of $M_{(1,1)}$. However, if the rank of $M_{(1,1)}$ is not full, then our analysis in \ref{lem:n-1,k-1} shows that $\vec{v}$ is not in the span of $M_{(1,1)}$, and further there are no restrictions on $\alpha$ (this is particular to symmetric matrices of full rank). Therefore, $\pi_{n,n}^{-1}(\Sym^{n-1,n-1})$ closed as it is cut out by the equation:
    \[
    \Bigg\{\begin{bmatrix}
        \alpha&\vec{v}\\
        \vec{v}&M_{(1,1)}
    \end{bmatrix}\in \Sym^{n,n}: \mathrm{det}\Bigg(\begin{bmatrix}
        0&\vec{v}\\
        \vec{v}&M_{(1,1)}
    \end{bmatrix}\Bigg)=0\Bigg\}
    \]

    This gives us: 
    \[
    M_{gm}(\Sym^{n,n})\cong M_{gm}(\pi_{n,n}^{-1}(\Sym^{n-1,n-1}))\oplus M_{gm}(\pi_{n,n}^{-1}(\Sym^{n-1,n-2}))
    \]

    Now, by the analysis in \ref{lem:n-1,k}, the projection 
    \[\pi_{n,n}^{-1}(\Sym^{n-1,n-1})\to \Sym^{n-1,n-1}\times\A^{n-1}\]
    sending 
    \[
    M\mapsto \begin{bmatrix}
        \vec{v}\\
        M_{(1,1)}
    \end{bmatrix}
    \]
    is a $\G_m$ bundle. so, it is the complement the zero section of an associated line bundle 
    $L\to \Sym^{n-1,n-1}\times\A^{n-1}$. One has the Gysin triangle:
    \[
    M_{gm}(\pi_{n,n}^{-1}(\Sym^{n-1,n-1}))\to M_{gm}(\Sym^{n-1,n-1})\to M_{gm}(\Sym^{n-1,n-1})(1)[2]\to M_{gm}(\pi_{n,n}^{-1}(\Sym^{n-1,n-1}))[1]
    \]
    where by 4.20 of \cite{Deglise} the map $M_{gm}(\Sym^{n, n})\to M_{gm}(\P\Sym^{n,n})(1)[2]$ is the Chern class map $(1_{\Sym^{n, n}})_*\boxtimes c_1(\gamma)$, where $c_1(\gamma):M_{gm}(\Sym^{n,n})\to \F(1)[2]$ is the first Chern class of $L$. We have (for any smooth variety)
\[
\mathrm{Hom}_{C^{gm}}(M_{gm}(\Sym^{n,n}),\F(1)[2])\cong \mathrm{Pic}(\Sym^{n,n})_\F\cong \mathrm{CH}^1(\Sym^{n, n})_\F
\]
which is zero as $\Sym^{n, n}$ is the complement of a closed subset of $\A^{{n\choose 2}}$. So, one can conclude that $c_1(\gamma)=0$, and so 
\[
M_{gm}(\pi_{n,n}^{-1}(\Sym^{n-1,n-1}))\cong M_{gm}(\Sym^{n-1,n-1})(1)[1]\oplus M_{gm}(\Sym^{n-1,n-1})
\]
concluding the result.\begin{footnote}    
    {Since $\mathrm{CH}^1(\Sym^{n, n})=0$, it follows that this splitting holds with \emph{any} coefficient ring $R$ if the base field is characteristic $0$, and for any ring $R$ such that $\mathrm{char}(\F)\in R$ is invertible if $\F$ is perfect of finite characteristic.} 
\end{footnote}
\end{proof}

In the case of symmetric matrices of full rank, the following can be said about their projectivization:
\begin{corollary}\label{projectivization}
    Let $\P\Sym^{n,n}$ denote the space of $n\times n$ projective symmetric matrices of rank $n$, and let the coefficient field $\F$ be of characteristic not equal to $n$. Then $M_{gm}(\P\Sym^{n,n})$ is split Tate, and 
    \[
M_{gm}(\Sym^{n,n})\cong M_{gm}(\P\Sym^{n,n})\otimes M_{gm}(\G_m)
\]
\end{corollary}
\begin{proof}
First, restrict the tautological bundle $\mathcal{O}_{\P^{n(n+1)/2}}(-1)$ to the open subvariety $\P\Sym^{n,n}$. The complement of the zero section of this bundle over $\P\Sym^{n,n}$ is exactly $\Sym^{n,n}$. So, we have the Gysin map:
\[
M_{gm}(\Sym^{n,n})\to M_{gm}(\P\Sym^{n,n})\to M_{gm}(\P\Sym^{n,n})(1)[2]\xrightarrow{+1}M_{gm}(\Sym^{n,n})[1]
\]
By the exact same argument as in \ref{inductive}, the map $M_{gm}(\P\Sym^{n,n})\to M_{gm}(\P\Sym^{n,n})(1)[2]$ vanishes because
\[\mathrm{CH}^1(\P\Sym^{n, n})_\F\cong \Z/n\Z\otimes\F=0
\]
So we have the splitting:
\[
M_{gm}(\Sym^{n,n})\cong M_{gm}(\P\Sym^{n,n})\oplus M_{gm}(\P\Sym^{n,n})(1)[1]
\]
Since $M_{gm}(\Sym^{n,n})$ is a finite direct sum of Tate motives, we can employ 2.73  of \cite{VoEilenberg} to say that $M_{gm}(\P\Sym^{n,n})$ is also a finite direct sum Tate motives as it is a direct summand. 
\end{proof}
\bibliographystyle{alpha}
\bibliography{refs}

\begin{thebibliography}{Dup24}

\bibitem[Blo86]{Bloch}
Spencer Bloch.
\newblock Algebraic cycles and higher k-theory.
\newblock {\em Advances in mathematics}, 61(3):267--304, 1986.

\bibitem[Bor14]{Borisov}
Lev Borisov.
\newblock Class of the affine line is a zero divisor in the grothendieck ring.
\newblock {\em arXiv preprint arXiv:1412.6194}, 2014.

\bibitem[Bro05]{Brosnan}
Patrick Brosnan.
\newblock On motivic decompositions arising from the method of
  bialynicki-birula.
\newblock {\em Inventiones mathematicae}, 161(1):91--111, 2005.

\bibitem[D{\'e}g08]{Deglise}
Fr{\'e}d{\'e}ric D{\'e}glise.
\newblock Around the gysin triangle. ii.
\newblock {\em Documenta Mathematica}, 13:613--675, 2008.

\bibitem[Dup24]{ClementDupont}
Cl{\'e}ment Dupont.
\newblock An introduction to mixed tate motives.
\newblock {\em arXiv preprint arXiv:2404.03770}, 2024.

\bibitem[Kel]{K}
Shane Kelly.
\newblock {\em Voevodsky motives and ldh-descent}.

\bibitem[Mac69]{M}
Jessie MacWilliams.
\newblock Orthogonal matrices over finite fields.
\newblock {\em The American Mathematical Monthly}, 76(2):152--164, 1969.

\bibitem[Sou85]{Soule}
Christophe Soul{\'e}.
\newblock Op{\'e}rations en k-th{\'e}orie alg{\'e}brique.
\newblock {\em Canadian Journal of Mathematics}, 37(3):488--550, 1985.

\bibitem[Tot14]{Totaro}
Burt Totaro.
\newblock Chow groups, chow cohomology, and linear varieties.
\newblock In {\em Forum of Mathematics, Sigma}, volume~2, page e17. Cambridge
  University Press, 2014.

\bibitem[Voe00]{Vo1}
Vladimir Voevodsky.
\newblock Triangulated categories of motives over a field.
\newblock {\em Cycles, transfers, and motivic homology theories}, 143:188--238,
  2000.

\bibitem[Voe10]{VoEilenberg}
Vladimir Voevodsky.
\newblock Motivic eilenberg-maclane spaces.
\newblock {\em Publications Math{\'e}matiques de l'IH{\'E}S}, 112:1--99, 2010.

\bibitem[Wil12]{Williams}
Ben Williams.
\newblock The motivic cohomology of stiefel varieties.
\newblock {\em Journal of K-theory}, 10(1):141--163, 2012.

\end{thebibliography}
\end{document}